\documentclass[10pt]{article}
\usepackage{graphicx, url}
\usepackage{lineno}
\usepackage{amsfonts}
\usepackage{amsmath}
\usepackage[linesnumbered,ruled,vlined]{algorithm2e}
\usepackage{algcompatible}
\usepackage{caption}

\modulolinenumbers[5]

\newcommand{\alined}{\emph{\textbf{aligned}}}

\newcommand{\cont}{\emph{\textbf{contained}}}
\newcommand{\scat}{\emph{\textbf{scattered}}}
\newcommand{\compat}{\emph{\textbf{compatible}}}

\newcommand{\alineds}{\emph{\textbf{aligned}} }
\newcommand{\nalineds}{\emph{\textbf{non-aligned}} }
\newcommand{\conts}{\emph{\textbf{contained}} }
\newcommand{\scats}{\emph{\textbf{scattered}} }
\newcommand{\compats}{\emph{\textbf{compatible}} }

\newcommand{\propone}{\textbf{$P_{1}$:} }
\newcommand{\proptwo}{\textbf{$P_{2}$:} }
\newcommand{\prboneref}{\textbf{$P_{1}$} }
\newcommand{\prbtworef}{\textbf{$P_{2}$} }

\newcommand{\prbtworefe}{\textbf{$P_{2}$}}

\newtheorem{thm}{Theorem}[section]
\newtheorem{cor}[thm]{Corollary}
 
 \newtheorem{prop}[thm]{Proposition}

\begin{document}

\begin{center}
{\LARGE Cache-oblivious Matrix Multiplication for Exact Factorisation}

\vspace{15pt} {\large Fatima K. Abu Salem\footnote{Corresponding author {\small E-mail:
fatima.abusalem@aub.edu.lb}}}

Computer Science Department, American University of Beirut, \\ P. O. Box 11-0236, Riad El Solh, Beirut 1107 2020, Lebanon

\vspace{15pt} {\large Mira Al Arab \footnote{{\small E-mail:
maa75@aub.edu.lb}}}

Computer Science Department, American University of Beirut, \\ P. O. Box 11-0236, Riad El Solh, Beirut 1107 2020, Lebanon
\end{center}

\begin{abstract}
We present a cache-oblivious adaptation of matrix multiplication to be incorporated in the parallel TU decomposition for rectangular matrices over finite fields, based on the Morton-hybrid space-filling curve representation. To realise this, we introduce the concepts of alignment and containment of sub-matrices under the Morton-hybrid layout. We redesign the decompositions within the recursive matrix multiplication to force the base case to avoid all jumps in address space, at the expense of extra recursive matrix multiplication (MM) calls. We show that the resulting cache oblivious adaptation has low span, and our experiments demonstrate that its sequential evaluation order demonstrates orders of magnitude improvement in run-time, despite the recursion overhead.
\end{abstract}

{\bf \it Keywords}: Locality of reference, Cache-oblivious Algorithms, Space-filling Curves, Morton-hybrid Layout, TU Decomposition, Finite Fields

\section{Introduction}

We present a cache-oblivious adaptation of matrix multiplication to be incorporated in the paralel TU decomposition for rectangular matrices over finite fields, based on the Morton-hybrid space-filling curve representation. Exact triangulisation of matrices is crucial for a large range of problems in Computer Algebra and Algorithmic Number Theory, where a basis of the solution set of the associated linear system is required. Our focal algorithm of reference is the TURBO algorithm of Dumas et al. \cite{Dumas} for exact LU decomposition. This algorithm recurses on rectangular and potentially singular matrices. TURBO significantly reduces the volume of communication on distributed systems, and retains optimal work and linear span. TURBO can also compute the rank in an exact manner. As benchmarked against some of the most efficient current exact elimination algorithms in the literature, TURBO incurs low synchronisation costs and reduces the communication cost featured in \cite{Ibarra80,Ibarra82} by a factor of one third when used with only one level of recursion on 4 processors. A significant part of TURBO consists of matrix factorisation, and so, adapting this kernel in a cache-oblivious fashion will ultimately contribute to a cache oblivious factorisation algorithm. That TURBO has low depth makes adapting its sequential version to the cache-oblivious model more telling. Particularly, nested parallel algorithms for which the natural sequential execution has low cache complexity will also attain good cache complexity on parallel machines with private or shared caches \cite{BlellAl.10}. 

At the base case of TURBO the sub-matrices reach a given threshold, and so one can take advantage of cache effects. To the best of our knowledge, no cache oblivious (or cache aware) algorithms for exact linear algebra exist in the literature. We pursue a cache oblivious adaptation using space-filling curves. TURBO requires index conversion routines from the space curve chosen and the cartesian order, due to the row and column permutations. In \cite{AAArx16a}, using a detailed analysis of the number of bit operations required for index conversion, and filtering the cost of lookup tables that represent the recursive decomposition of the Hilbert curve, we have shown that the Morton-hybrid order incurs the least cost for index conversion routines as compared to the Hilbert, Peano, or Morton orders. The Morton order is the recursive $Z$-shaped space filling curve (Fig. \ref{fig:fmorton_refinement}). The Morton-hybrid order stops decomposing when the submatrices attain a threshold dimension $T' \times T'$ \cite{Wise:MaskedIntegers}. At such a level, say when the submatrix fits in cache, the overhead for maintaining the curve representation outweighs the reduction in cache complexity. In reference to the literature cited in this manuscript around the Morton-order and its hybrid, this curve representation improves significantly on the temporal locality of various matrix algorithms such as naive multiplication, LU decomposition, and QR factorisation. 

In this work, we introduce the concepts of alignment and containment of sub-matrices under the Morton-hybrid layout, and develop the full details of the MM algorithm by which it observes the alignment and containment of sub-matrices invariably across the matrix factorisation recursive steps. We do this by redesigning the decompositions within the recursive MM to force the base case to avoid all jumps in address space, at the expense of extra recursive MM calls. We show that the resulting cache oblivious adaptation retains optimal work and critical path length as default MM and thus is highly parallel. Our experiments confirm that the recursion overhead in the Morton-hybrid MM is negligible and leads to significant reduction in run-time thanks to its improved temporal locality. 

Before proceeding, we begin with brief description of the TU algorithm.  Consider a rectangular matrix $A$ over a field $\mathbb{F}$, where $A$ may be singular. $A$ is triangulated into the product of two matrices $T$ and $U$, such that $A=T \cdot U$, where $U$ is a upper triangular matrix, and $T$ is with some ``$T$'' patterns. This is done in a series of recursive steps on rectangular and potentially singular matrices, relaxing the condition for generating a strictly lower triangular matrix: (1) Recursive TU decomposition in SE, SW, NE, and NW (2) Virtual row and column permutations needed to re-order the blocks to yield the matrix $U$. For brevity and because of lack of space, we omit the full details of the algorithm and refer the reader to \cite{Dumas} for a full account on TURBO.


\section{Non-Aligned Rectangular Sub-matrix Multiplication Within The Recursion}
\label{sec:rec_mult}

Consider Morton-hybrid matrices $A$, $B$, and $C$ and let $S_{A}$, $S_{B}$, and $S_{C}$ be random sub-matrices of $A$, $B$, and $C$ respectively, for which one has to compute $S_{A} = S_{B} \cdot S_{C}$. This is a typical scenario encountered during the TU decomposition. To illustrate further, consider Fig. \ref{fig:submm_example}. Each integer appearing in the matrices in that figure represents the corresponding Morton-hybrid index of the element occupying it. The sub-matrices on which the multiplication is performed do not begin at the first entry of a Morton-hybrid sub-matrix, hence the concept of an \alineds versus \nalineds Morton-hybrid sub-matrix.

\begin{def}
\label{def:morton_aligned}An \textbf{aligned} sub-matrix is a
$2^{a}\cdot T \times 2^{b}\cdot T$ sub-matrix of a Morton-hybrid
matrix that begins at the first entry of a row-major sub-matrix. A
\textbf{non-aligned} sub-matrix is a sub-matrix of a Morton-hybrid
matrix that does not satisfy this condition.
\end{def}

\begin{cor}
\label{cor:alined} The Cartesian index of the first entry of an
\textbf{aligned} sub-matrix is of the form $(k_{1}\cdot T,
k_{2}\cdot T)$, for any positive integers $k_{1}$ and $k_{2}$.
\end{cor}

{\it Proof}: By its definition, an \alineds sub-matrix
$A_{M}$ of a Morton-hybrid matrix $M$ starts at the first entry of
some row-major sub-matrix $S_{M}$ of $M$. Since the row-major
sub-matrix $S_{M}$ is of dimensions $T' \times T'$, the Cartesian
index of the first entry of $S_{M}$ is given by $(k_{1} \cdot T,
k_{2} \cdot T)$, for some positive integers $k_{1}$ and $k_{2}$.
By its definition, the \alineds sub-matrix $A_{M}$ begins
at an element of Cartesian index $(k_{1}\cdot T, k_{2}\cdot T)$.

\begin{cor}
\label{cor:alined_rowmajor}If an \textbf{aligned} sub-matrix is $T' \times T'$, then it is row-major.
\end{cor}

{\it Proof}: 
Let $A_{M}$ be a $T' \times T'$ \alineds sub-matrix of a Morton-hybrid
matrix $M$. From the definition of an aligned matrix, we know that $A_{M}$
begins at the first entry of a row-major sub-matrix $S_{M}$ of $M$.
According to the Morton-hybrid layout, all row-major sub-matrices of
$M$, including $S_{M}$ are $T' \times T'$. Since $A_{M}$ is also $T'
\times T'$, then $A_{M}$ must be $S_{M}$ and hence is row-major.

An example of a \nalineds sub-matrix of a Morton-hybrid matrix with
$T'=4$ is shown in red in Fig. \ref{fig:non_aligned}. An \alineds
sub-matrix is shown in green.

\noindent Next, we relate the lack of alignment of sub-matrices to
the recursive accessing of these sub-matrices and discuss the
implicated problems.

\subsection{Non-Aligned Sub-Matrices and loss of locality}
\label{sec:non_aligned_problems} 
\begin{def}
\label{def:contained}A sub-matrix $S_{M}$ of a Morton-hybrid matrix
$M$ is said to be \textbf{contained} if $S_{M}$ lies completely
within a sub-matrix of $M$ ordered in a row-major fashion.
Otherwise, we say that $S_{M}$ is \textbf{scattered}.
\end{def}

\begin{prop}\label{prps:alined_contained}
Let $A_{M}$ be an \alineds sub-matrix of a Morton-hybrid matrix $M$.
The sub-matrix at the base case of the recursive division, down
until $T' \times T'$ sub-matrices, of $A_{M}$ is a $T' \times T'$
row-major sub-matrix of $M$.
\end{prop}

{\it Proof}: 
First, we claim the recursive division of $A_{M}$ gives 4 \alineds sub-matrices. From the definition of \alineds sub-matrices, $A_{M}$ has size $2^{a}\cdot T \times
2^{b}\cdot T$. So, the division of each of the dimensions of $A_{M}$
by 2 results in four quadrants $NW$, $NE$, $SW$, and $SE$ of $A_{M}$
of size $2^{(a-1)}\cdot T \times 2^{(b-1)}\cdot T$ each. Thus these
quadrants satisfy the size condition from
the definition of aligned matrices. Note that once any of the dimensions
reaches size $T'$ it is no longer divided, and the recursive division
proceeds on the other dimension until that too becomes $T'$. It is
the size condition of this same definition that leads to $T'
\times T'$ sub-matrices at the base case of recursive division of the
\alineds sub-matrices decomposed from $A_{M}$.

Now, recall, from Cor. \ref{cor:alined},
that the start index of $A_{M}$ is of the form $(k_{1}\cdot T,
k_{2}\cdot T)$. Then, the start indices of the sub-matrices
resulting from the sub-division of $A_{M}$ are $(k_{1}\cdot T,
k_{2}\cdot T)$, $(k_{1}\cdot T, (k_{2} + 2^{(b-1)})\cdot T)$, $(
(k_{1}+ 2^{(a-1)})\cdot T, k_{2}\cdot T)$, and $(k_{1}+
2^{(a-1)})\cdot T, (k_{2} + 2^{(b-1)})\cdot T)$ for the $NW$, $NE$,
$SW$, and $SE$ quadrants of $A_{M}$ respectively. Thus the start
indices of these quadrants satisfy the start index condition from
the definition of aligned matrices. Combining, by those two claims, the four quadrants
resulting from the sub-division of any \alineds matrix $A_{M}$ are
\alineds: they satisfy both conditions from
the definition of aligned matrices. 

Second, we show that the \alineds sub-matrices at the base case are
row-major sub-matrices of $M$. If the recursive division continues
till $T' \times T'$ sub-matrices, we get $T' \times T'$ \alineds
sub-matrices. From Cor. \ref{cor:alined_rowmajor}, we know that these
sub-matrices are row-major sub-matrices of $M$. This concludes the
proof.

\begin{cor}
\label{cor:alined_contained} Any sub-matrix of the $T' \times T'$
sub-matrix reached at the base case of the recursive division of an
\textbf{aligned} sub-matrix is \textbf{contained}.
\end{cor}

{\it Proof}: 
According to Prop. \ref{prps:alined_contained}, the $T' \times T'$
sub-matrix at the base case of the recursive division of an \alineds
sub-matrix $A_{M}$ of a Morton-hybrid matrix $M$ is in row-major
layout. Hence, any sub-matrix of this $T' \times T'$ base
case sub-matrix lies entirely within a row-major sub-matrix of $M$
and is therefore \cont.

In Fig. \ref{fig:non_aligned_basecase}, $C_{M}$ is one of the sub-matrices at the base case of the recursive
division of the \alineds sub-matrix $A_{M}$ and is a row-major sub-matrix. Any sub-matrix of
$C_{M}$ is \cont. When \nalineds sub-matrices are recursively divided, the sub-matrix
at the base case may not consist entirely of a row-major sub-matrix
of the Morton-hybrid matrix. It may be \scats across more than one
row-major sub-matrix. For example, in Fig. \ref{fig:non_aligned_basecase}, the sub-matrix $S_{M}$ is a sub-matrix at the base case of the recursion for the
\nalineds sub-matrix $N_{M}$ in red. $S_{M}$ spans four row-major ordered sub-matrices - hence, it is \scat. We know that the elements of the sub-matrices at the base case are
to be traversed in a row-major or column-major order, as required
for the base case of MM. With such traversal imposed, a \scats
sub-matrix suffers from two issues:
\begin{enumerate}
\item{
\propone \textit{Elements of a \textbf{scattered} sub-matrix are
not sufficiently close in memory to maintain good spatial locality
when traversed in a row/column-major fashion. This results in
worse memory performance than for \textbf{contained} sub-matrices.}}
\item{
\proptwo \textit{Morton-hybrid encoding is required for accessing
each element within a \textbf{scattered} sub-matrix (thus incurring
extra computation overhead compared to row-major offset
calculation).}}
\end{enumerate}
\begin{prop}
The loss in locality defined by \prboneref and \prbtworef apply
for \scats sub-matrices but not \conts sub-matrices.
\end{prop}
{\it Proof}: 
We first consider \prboneref. Recall that the traversal of entries at the base case of
the recursion is done in two orders: row-major and column-major. For
\conts sub-matrices, when consecutively accessing any two entries in
any of these two orders, the minimum jump in address space is 1 and
the maximum is $T'$ as all entries lie within one row-major
sub-matrix of the Morton-hybrid matrix. A \scats sub-matrix spans
more than one row-major sub-matrix of the matrix. These row-major
sub-matrices are not necessarily consecutive in memory and
traversing, in a row-major or column-major fashion, the \scats
sub-matrix that spans these row-major sub-matrices results in jumps
in address space. When consecutively accessing any two entries of a
\scats sub-matrix, the minimum jump in address space is 1 if the two
entries being accessed consecutively belong to the same row-major
sub-matrix and the maximum is $k \cdot T^{2} + T - 1$ for some
positive integer $k$, if the two entries belong to different
row-major sub-matrices.

We now consider \prbtworefe. Because the base case sub-matrix of an
\alineds sub-matrix is part of a row-major ordered sub-matrix,
offset calculation for the elements at the base case is fast:
traditional row-major offset calculation is used. Index $z$ of an
element at offset $(i,j)$ from the start index $\sigma$ of the
sub-matrix at the base case is given by $z = \sigma + i*T + j$,
since the sub-matrix satisfies a row-major ordering with row length = $T'$. This can be seen
for the \conts sub-matrix $C_{M}$ shown in
Fig. \ref{fig:non_aligned_basecase} where $\sigma = 112$. As for a
\nalineds sub-matrix, accessing any element ($i,j$) in any of the
base case sub-matrices requires that the corresponding Morton-hybrid
index be calculated. This incurs extra calculation overhead as the
encoding of the Morton-hybrid index is more costly than calculating
an offset within a row-major ordered sub-matrix.

\subsection{Modified Non-Aligned Sub-Matrix Multiplication}
\label{sec:my_approach} 
\begin{figure}[!tbp]
  \centering
     {\includegraphics[width=3in]{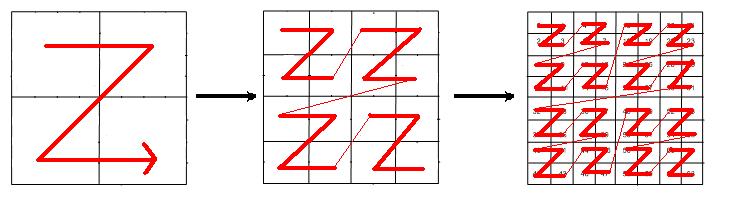}}
  \caption{Generation of Morton order}
  \label{fig:morton_refinement}
  \hfill
  \begin{minipage}[b]{0.4\textwidth}
     {\includegraphics[width=2.3in]{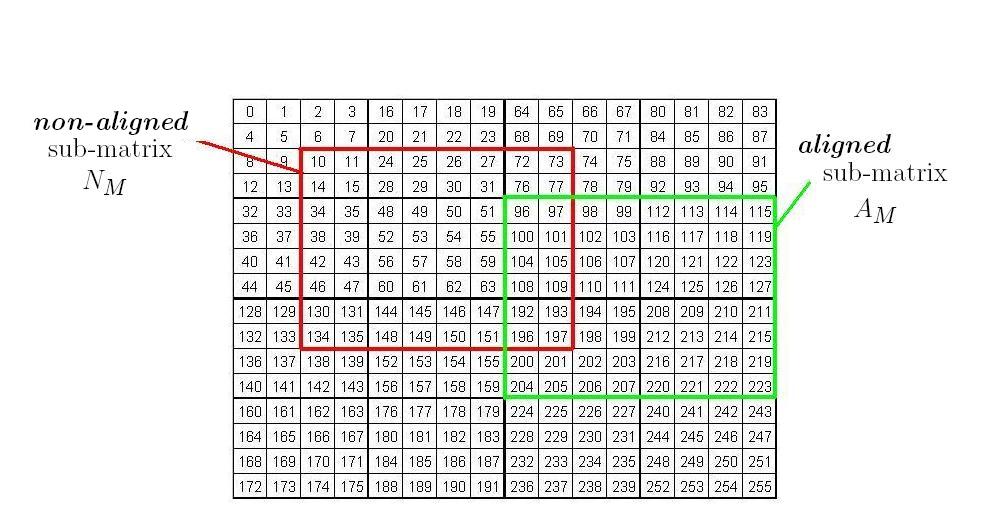}}
  \caption{Aligned and non-aligned sub-matrices}
  \label{fig:non_aligned}     
  \end{minipage}
\hfill
  \begin{minipage}[b]{0.4\textwidth}
{\includegraphics[width=3in]{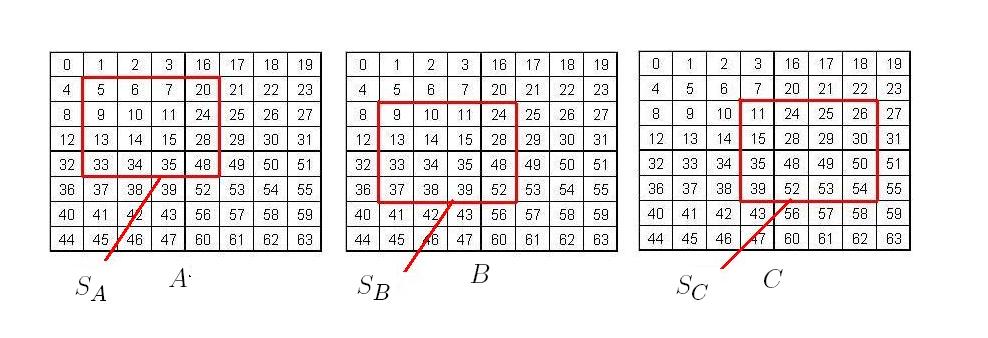}}
  \caption{Sub-matrix multiplication example I}
  \label{fig:submm_example}
  \end{minipage}
  \hfill
  \begin{minipage}[b]{0.4\textwidth}
     {\includegraphics[width=2.5in]{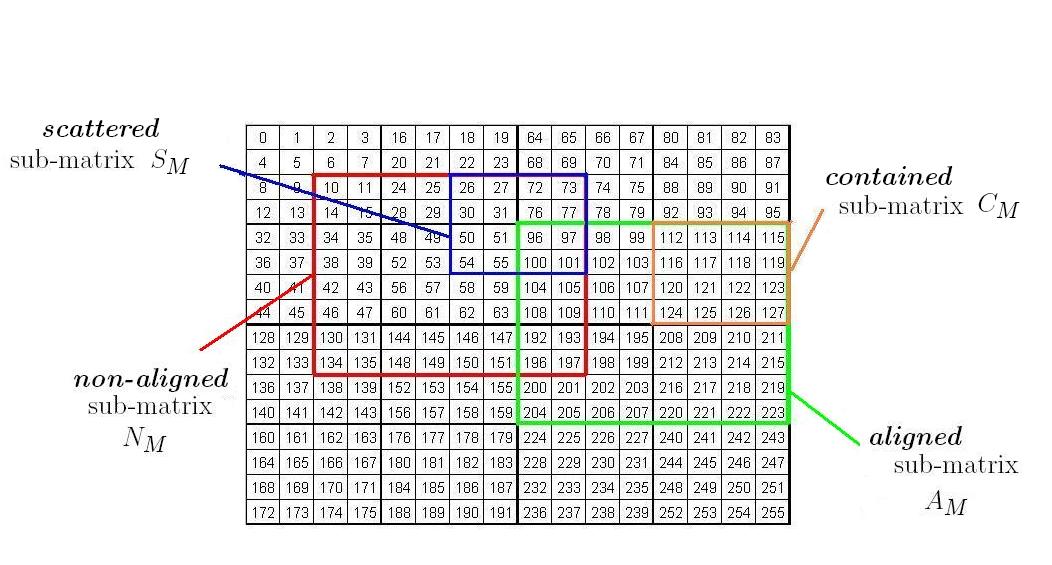}}
  \caption{Base case sub-matrix of non-aligned sub-matrix}
  \label{fig:non_aligned_basecase}
  \end{minipage}
  \hfill
  \begin{minipage}[b]{0.4\textwidth}
     {\includegraphics[width=2in]{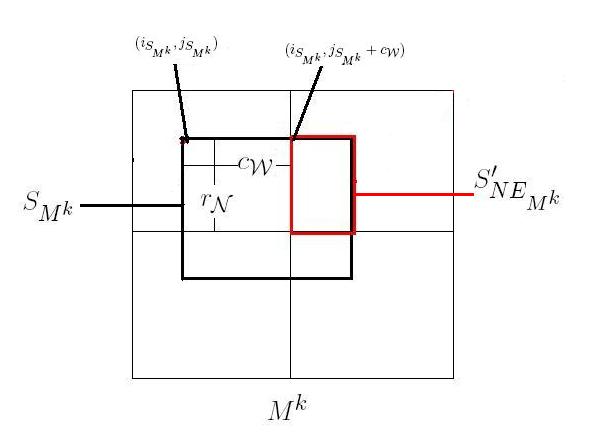}}
  \caption{$S'_{NE_{M^{k}}}$}
  \label{fig:submm_example_cp}
  \end{minipage}
  \hfill
     {\includegraphics[width=3in]{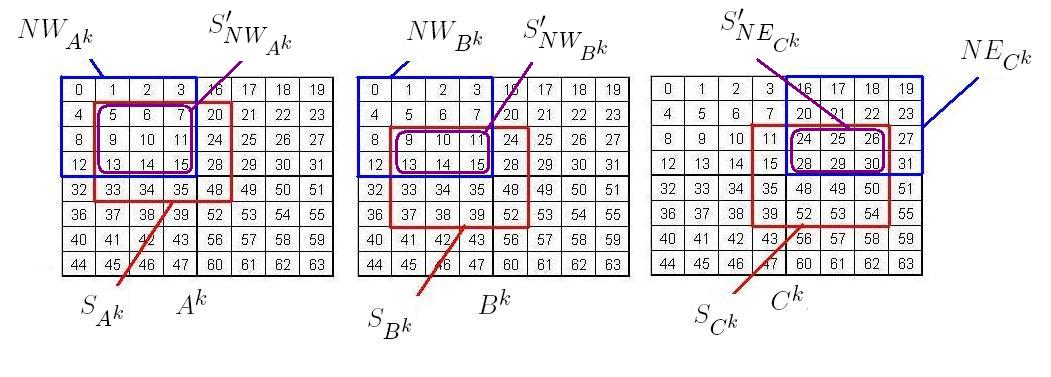}}
  \caption{Sub-matrix multiplication example II}
  \label{submatrix_mult_example_mh1}
  \begin{minipage}[b]{0.4\textwidth}
     {\includegraphics[width=3in]{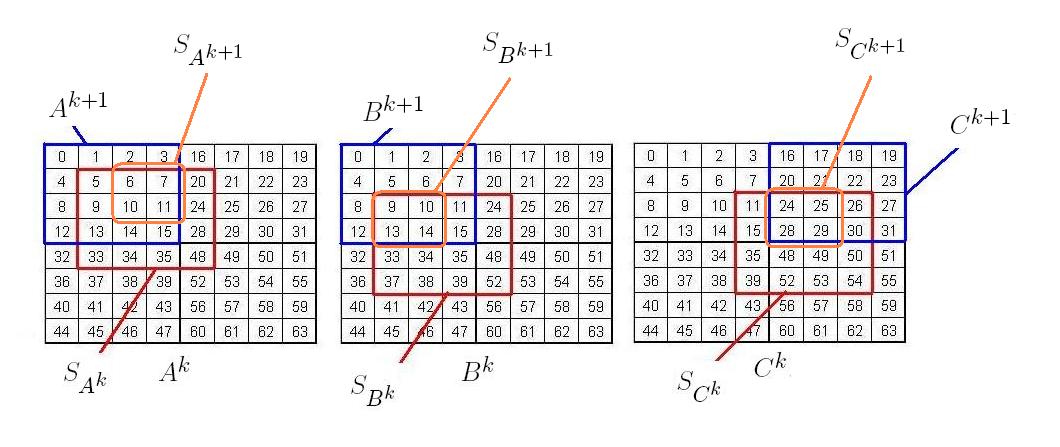}}
  \caption{Example of Modified Recursive Sub-matrix Multiplication}
  \label{fig:example_sprime} 
  \end{minipage}
\end{figure}

We aim to improve the sub-matrix
multiplication procedure by addressing issues \prboneref and
\prbtworef. In this section, we describe a recursive sub-matrix
multiplication algorithm which ensures that the sub-matrices at the
base case of the recursion are \conts in a row-major ordered
sub-matrix of the original matrix. By doing this, we reduce the
range of addresses of the elements within the sub-matrices at the
base case as well as the number of jumps in address space done at
the base case, and we eliminate the need for Morton-hybrid encoding
at the base case. To ensure efficiency that the sub-matrix at the
base case of MM is \cont, by Prop. \ref{prps:alined_contained}, the
recursive division within the algorithm must start on \alineds
matrices. Recall the random matrices $A$, $B$, and $C$ in Morton-hybrid order and of dimensions $2^{m}
\times 2^{m}$, and $S_{A}$, $S_{B}$, and $S_{C}$ the random
sub-matrices of $A$, $B$, and $C$ respectively (Fig. \ref{fig:submm_example}). We wish to
perform the multiplication $S_{A} = S_{B} \cdot S_{C}$ efficiently. We can
recursively divide $S_{A}$, $S_{B}$, and $S_{C}$, as in
the default MM algorithm, which
may result in \scats sub-matrices at the base case since $S_{A}$,
$S_{B}$, and $S_{C}$ may not be \alined. Instead, we will
recursively divide $A$, $B$, and $C$ and address only the relevant
sub-matrix multiplications that ought to be done to produce $S_{A} =
S_{B} \cdot S_{C}$. As $A$, $B$, and $C$ are \alined, recursively
dividing them will enforce row-major sub-matrices at the base case
from which we extract the relevant parts to produce $S_{A}$.

Let $k$ be a superscript denoting a recursive step of the proposed MM algorithm. Also, let $t$, $u$, and $v$, denote subscripts in $\{0,1,2,3\}$ (of sub-matrices of $A$, $B$, and $C$ respectively), indicating a specific quadrant following the Morton (Z-order): $NW = 0$, $NE = 1$, $SW = 2$, and $SW = 3$. For $k=0$, $A^{0}_{0} = A$, $B^{0}_{0} = B$, and $C^{0}_{0} = C$. Denote by $S_{A^{k}_{t}}$, $S_{B^{k}_{u}}$, and and $S_{C^{k}_{v}}$ the respective sub-matrices of $A^{k}_{t}$, $B^{k}_{u}$, and $C^{k}_{v}$ being multiplied as part of the overall multiplication $S_{A} = S_{B} \cdot S_{C}$. As such the initial problem is to produce $S_{A^{0}_{0}} = S_{B^{0}_{0}} \cdot S_{C^{0}_{0}}$. For this, we first produce the quadrants $A^{1}_{t'}$ of $A^{0}_{0}$, such that $A^{1}_{t'} \in \{NW_{A^{0}_{0}},NE_{A^{0}_{0}},SW_{A^{0}_{0}},SE_{A^{0}_{0}}\}$ for $t' \in \{0,1,2,3\}$. We do the same for $B^{0}_{0}$ and $C^{0}_{0}$ producing $B^{1}_{u'}$ and $C^{1}_{v'}$ respectively for $u',v' \in \{0,1,2,3\}$. For each $A^{1}_{t'}$, we produce the sub-matrix$S'_{A^{1}_{t'}}$, defined as the part of $S_{A^{0}_{0}}$ that lies in $A^{1}_{t'}$. Similarly, we produce $S'_{B^{1}_{u'}}$, and $S'_{C^{1}_{v'}}$. Note that $S_{A^{0}_{0}}$ is the two-dimensional concatenation of $\{S'_{A^{1}_{t'}}\}$ for $t' \in \{0,1,2,3\}$ and hence to calculate $S_{A^{0}_{0}}$ we need to calculate $S'_{A^{1}_{t'}}$ for $t' \in \{0,1,2,3\}$. To do this, we need to consider all combinations $\Gamma$ of the form $\Gamma_{t',u',v'} = \{A^{1}_{t'}, B^{1}_{u'}, C^{1}_{v'}\}$ necessary to produce $S_{A^{0}_{0}}$, as will be justified below. Now, when considering a combination $\Gamma_{t',u',v'} = \{A^{1}_{t'}, B^{1}_{u'}, C^{1}_{v'}\}$, if the sub-matrices $S'_{A^{1}_{t'}}$, $S'_{B^{1}_{u'}}$, and $S'_{C^{1}_{v'}}$ are \compats for multiplication, i.e. the multiplication $S'_{A^{1}_{t'}} += S'_{B^{1}_{u'}} \cdot S'_{C^{1}_{v'}}$ is part of the overall multiplication $S_{A^{0}_{0}} += S_{B^{0}_{0}} \cdot S_{C^{0}_{0}}$, then a recursive call is made on $S'_{A^{1}_{t'}}$, $S'_{B^{1}_{u'}}$, and $S'_{C^{1}_{v'}}$. Else, if $S'_{A^{1}_{t'}}$, $S'_{B^{1}_{u'}}$, and $S'_{C^{1}_{v'}}$ are not \compat, we extract \emph{\textbf{compatible}} parts of these sub-matrices and we label them as $S_{A^{1}_{t'}}$, $S_{B^{1}_{u'}}$, and $S_{C^{1}_{v'}}$ on which the multiplication proceeds recursively. After doing this for all combinations $\Gamma_{t',u',v'}$ for $t',u',v' \in \{0,1,2,3\}$, we would have calculated $S_{A^{0}_{0}}$. 

We now describe the general $k$'th recursive step of Morton-hybrid MM, which consists of a round of four substeps. For simplicity, we drop the subscripts $t$, $u$ and $v$ of $A^{k}_{t}$, $B^{k}_{u}$, and $C^{k}_{v}$, and we use $M$ to denote any of the matrices $A$, $B$, or $C$, Each \alineds $M^{k}$ is identified by two values:
\begin{description}
\item[$\alpha_{M^{k}}$:] the Morton-hybrid index of the first element in the \alineds matrix $M^{k}$
\item[$\lambda_{M^{k}}$:] the number of elements in the \alineds sub-matrix $M^{k}$
\end{description}
We are also given the sub-matrices $S_{A^{k}}$ of $A^{k}$,
$S_{B^{k}}$ of $B^{k}$, and $S_{C^{k}}$ of $C^{k}$ on which we wish
to perform the multiplication. Each of the sub-matrices $S_{M^{k}}$
is identified by the following:
\begin{description}
\item[$\sigma_{S_{M^{k}}}$:] the Morton-hybrid index of the first entry of $S_{M^{k}}$
\item[$r_{S_{M^{k}}}$:] the number of rows of $S_{M^{k}}$
\item[$c_{S_{M^{k}}}$:] the number of columns of $S_{M^{k}}$
\end{description}
We do not use the 4-tuple $(M,\sigma,r,c)$ to identify the \alineds
sub-matrices $M^{k}$ because the 3-tuple $(M,\alpha,\lambda)$
simplifies the computations for identification of the quadrants of
$M^{k}$ and incorporates the information from the 4-tuple where
$\alpha = \sigma$ and $\lambda = r \times c$.

\noindent {\bf Step 1:} In this step, we need to identify all four \alineds quadrants $M^{k+1}_{t}$, $t \in
\{0,1,2,3\}$, of the \alineds $M^{k}$, for $k$ not reaching the base
case, to proceed with the recursive multiplication algorithm. The
index $t$ is dropped from $M_{k}$ for simplicity. To do this, we
identify the start index $\alpha_{M^{k+1}_{t}}$ and size
$\lambda_{M^{k+1}_{t}}$ of each quadrant $M^{k+1}_{t}$ of $M^{k}$.
Because $M^{k}$ is divided into four quadrants of equal size, the
number of elements $\lambda_{M^{k+1}_{t}}$ in any quadrant is given
by $\lambda^{\prime} = \lambda_{M^{k+1}_{t}} = \lambda_{M^{k}} / 4$. 
Recall, that in the Morton-hybrid order, the quadrants not reaching
the base case are stored according to the Morton layout. For the Morton order, the quadrants of $M^{k}$ are laid out in the order
$NW_{M^{k}}$, $NE_{M^{k}}$, $SW_{M^{k}}$, then $SE_{M^{k}}$, and
hence
\begin{itemize}
\item{$\alpha_{NW_{M^{k}}} = \alpha_{M^{k}}$}
\item{$\alpha_{NE_{M^{k}}} = \alpha_{M^{k}} + \lambda^{\prime}$}
\item{$\alpha_{SW_{M^{k}}} = \alpha_{M^{k}} + 2(\lambda^{\prime})$}
\item{$\alpha_{SE_{M^{k}}} = \alpha_{M^{k}} + 3(\lambda^{\prime})$}
\end{itemize}
The sub-matrix $S_{M^{k}}$ may not lie entirely within one quadrant
of $M^{k}$ and hence all quadrants $M^{k+1}_{t}$ of $M^{k}$ which
contain part of $S_{M^{k}}$ must be considered, which is the case in
the example from Fig. \ref{fig:submm_example_cp} as $S_{A^{k}}$,
$S_{B^{k}}$, and $S_{C^{k}}$ touch on all four quadrants of $A^{k}$,
$B^{k}$, and $C^{k}$ respectively. Given $S_{M^{k}}$, we must now
identify, for each $M^{k+1}_{t}$, the part of $S_{M^{k}}$ that lies
within $M^{k+1}_{t}$. We denote this sub-matrix by
$S'_{M^{k+1}_{t}}$. The method to identify $S'_{M^{k+1}_{t}}$ now follows.

\noindent {\bf Step 2:} Recall that we are given $M^{k}$ and $S_{M^{k}}$ as input into the
recursion. As $S_{M^{k}}$ may not lie entirely within
one quadrant of $M^{k}$, it is \scat, and we must identify the parts
of $S_{M^{k}}$ which lie in $M^{k+1}$ denoted by
$S^{\prime}_{M^{k+1}}$. We have identified the quadrants $M^{k+1}_{t}$, for $t
\in \{0,1,2,3\}$, and now we will identify the part of $S_{M^{k}}$
that lies within each $M^{k+1}_{t}$, denoted by $S'_{M^{k+1}_{t}}$.
Then $S_{M^{k}}$ is the two-dimensional concatenation of
$\{S'_{M^{k+1}_{t}}\}$ for $t \in \{0,1,2,3\}$. Here we drop the
index $t$ for simplicity. To identify $S'_{M^{k+1}}$, we need to
identify its start index $\sigma_{S'_{M^{k+1}}}$ and dimensions $r_{S'_{M^{k+1}}} \times c_{S'_{M^{k+1}}}$. To do this, the following intermediate values are
needed. For simplicity, the indices of the intermediate values
denoting dependence on $M^{k}$ are omitted.

\begin{description}
\item[$r_{\mathcal{N}}$:] The number of rows of $S_{M^{k}}$ in $\mathcal{N}$, the northern half of $M^{k}$.
\item[$c_{\mathcal{W}}$:] The number of columns of $S_{M^{k}}$ in $\mathcal{W}$, the western half of $M^{k}$.
\item[$r_{\mathcal{S}}$:] The number of rows of $S_{M^{k}}$ in $\mathcal{S}$, the
southern half of $M^{k}$. 
\item[$c_{\mathcal{E}}$:] The number of columns of $S_{M^{k}}$ in $\mathcal{E}$, the
eastern half of $M^{k}$. 
\item[$e_{NW_{M^{k}}}$:] The Morton-hybrid index of the last entry of
$NW_{M^{k}}$. Similarly for $NE_{M^{k}}$ , $SW_{M^{k}}$, and
$SE_{M^{k}}$.
\item[$encode(i,j)$:] Given an entry $e$ of Cartesian index
$(i,j)$, $encode(i,j)$ returns the Morton-hybrid index of $e$
\item[$extract\_i(z)$:] Given an entry $e$ of Morton-hybrid index $z$, $extract\_i(z)$ returns the
coordinate $i$ of the Cartesian index $(i,j)$ of $e$
\item[$extract\_j(z)$:] Given an entry $e$ of Morton-hybrid index $z$, $extract\_j(z)$ returns the
coordinate $j$ of the Cartesian index $(i,j)$ of $e$
\end{description}
The identification of $S'_{M^{k+1}}$ is done as follows:
\begin{itemize}
\item{For $NE_{M^{k}}$, calculate $e_{NE_{M^{k}}}$ as follows:
\begin{equation*}
e_{NE_{M^{k}}} = \alpha_{NE_{M^{k}}} + \lambda' - 1 = \alpha_{M^{k}}
+ 2\cdot\lambda^{\prime} - 1
\end{equation*}
and, for $SW_{M^{k}}$, $e_{SW_{M^{k}}}$ use
\begin{equation*}
e_{SW_{M^{k}}} = \alpha_{SW_{M^{k}}} + \lambda' - 1 = \alpha_{M^{k}}
+ 3\cdot\lambda^{\prime} - 1.
\end{equation*} }
\item{Find $r_{\mathcal{N}}$ using $r_{\mathcal{N}} = extract\_i(e_{NE_{M^{k}}}) -
extract\_i(\sigma_{S_{M^{k}}}) + 1$, i.e. $r_{\mathcal{N}}$ is the difference between the row indices of
the last entry of $NE_{M^{k}}$ and the first entry of $S_{M^{k}}$
and represents the number of rows of $S_{M^{k}}$ in the northern
half of $M^{k}$. Similarly, we find
\begin{equation*}
c_{\mathcal{W}} = extract\_j(e_{SW_{M^{k}}}) -
extract\_j(\alpha_{S_{M^{k}}}) + 1,
\end{equation*} the number of columns of $S_{M^{k}}$ in the western half of $M^{k}$. Note that if $r_{\mathcal{N}} <= 0$ then
no part of $S_{M^{k}}$ lies in the north half of $M^{k}$ and if
$c_{\mathcal{W}} <= 0$ then no part of $S_{M^{k}}$ lies in the west
half of $M^{k}$. After finding $r_{\mathcal{N}}$ and
$c_{\mathcal{W}}$, we can find $r_{\mathcal{S}}$ and
$c_{\mathcal{E}}$ using $r_{\mathcal{S}} = r_{S_{M^{k}}} - r_{\mathcal{N}}$ and $c_{\mathcal{E}} = c_{S_{M^{k}}} -
c_{\mathcal{W}}$, which are the remaining rows and columns of $S_{M^{k}}$ respectively}
\item{ So far, we have found the number of rows in the northern and
southern halves of $M^{k}$ and the number of columns in the western
and eastern halves of $M^{k}$ and we want to identify
$S'_{NW_{M^{k}}}$, $S'_{NE_{M^{k}}}$, $S'_{SW_{M^{k}}}$ and
$S'_{SW_{M^{k}}}$ for each $M^{k} \in \{A^{k},B^{k},C^{k}\}$. Recall
that we are able to identify a sub-matrix by
a 4-tuple $(M,\sigma,r,c)$, where $\sigma$ is the Morton-hybrid
index of the first entry of the sub-matrix and $r$ and $c$ are its
row and column dimensions respectively. Let
$(i_{S_{M^{k}}},j_{S_{M^{k}}})$ denote the Cartesian index of the
first entry of $S_{M^{k}}$ found using $i_{S_{M^{k}}} = extract\_i(\sigma_{S_{M^{k}}})$ and
$j_{S_{M^{k}}} = extract\_j(\sigma_{S_{M^{k}}})$. 
We now identify $S'_{NW_{M^{k}}}$, $S'_{NE_{M^{k}}}$,
$S'_{SW_{M^{k}}}$ and $S'_{SW_{M^{k}}}$ according to the following
cases:

\begin{enumerate}
\item{\label{id_nw}For $S'_{NW_{M^{k}}}$, $\sigma_{S'_{NW_{M^{k}}}} = \sigma_{S_{M^{k}}}$, $r_{S'_{NW_{M^{k}}}} = r_{\mathcal{N}}$, and $c_{S'_{NW_{M^{k}}}} = c_{\mathcal{W}}$.}
\item{\label{id_ne}For $S'_{NE_{M^{k}}}$, $\sigma_{S'_{NE_{M^{k}}}} = encode(i_{S_{M^{k}}}, j_{S_{M^{k}}}+c_{\mathcal{W}})$, $r_{S'_{NE_{M^{k}}}} = r_{N}$, and $c_{S'_{NE_{M^{k}}}} = c_{\mathcal{E}}$.}
\item{\label{id_sw}For $S'_{SW_{M^{k}}}$, $\sigma_{S'_{SW_{M^{k}}}} = encode(i_{S_{M^{k}}}+r_{\mathcal{N}}, j_{S_{M^{k}}})$, $r_{S'_{SW_{M^{k}}}} = r_{\mathcal{S}}$, and $c_{S'_{SW_{M^{k}}}} = c_{\mathcal{W}}$.}
\item{\label{id_se}For $S'_{SE_{M^{k}}}$, $\sigma_{S'_{SE_{M^{k}}}} = encode(i_{S_{M^{k}}}+r_{\mathcal{N}}, j_{S_{M^{k}}}+c_{\mathcal{W}})$, $r_{S'_{SE_{M^{k}}}} = r_{\mathcal{S}}$, and $c_{S'_{SE_{M^{k}}}} = c_{\mathcal{E}}$.}
\end{enumerate}
}
\end{itemize}

To justify these cases we will explain how we arrived at case
\ref{id_ne} for example where we identify the start index and
dimensions of $S'_{NE_{M^{k}}}$ as shown in
Fig. \ref{fig:example_sprime}. The rest follow similarly. Recall that
$\sigma_{S_{M^{k}}}$ denotes the Morton-hybrid index of the first
element of $S_{M^{k}}$, and that $(i_{S_{M^{k}}},j_{S_{M^{k}}})$ is
the corresponding Cartesian index. The index
$(i_{S_{M^{k}}},j_{S_{M^{k}}}+c_{\mathcal{W}})$ is the Cartesian
index of the first element in $S'_{NE_{M^{k}}}$. The corresponding
Morton-hybrid index $\sigma_{S'_{NE_{M^{k}}}}$ can be found using
the function $encode(i_{S_{M^{k}}}, j_{S_{M^{k}}}+c_{\mathcal{W}})$.
The dimensions of $S'_{NE_{M^{k}}}$ are $r_{\mathcal{N}} \times
c_{\mathcal{E}}$.

Note that for each $S'_{M^{k+1}}$, the Cartesian index of the start
entry of Morton-hybrid index $\sigma_{S'_{M^{k+1}}}$ is given by
$(i_{S_{M^{k}}}+\varphi_{r_{M}} ,j_{S_{M^{k}}}+\varphi_{c_{M}})$,
for $\varphi_{r_{M}} \in \{0,r_{\mathcal{N}}\}$ and $\varphi_{c_{M}}
\in \{0,c_{\mathcal{W}}\}$. 

\noindent {\bf Step 3:} By now we have decomposed each $M^{k}$ into quadrants, and we have identified, for each quadrant $M^{k+1}_{t}$, the part of $S_{M^{k}}$ within that quadrant denoted by $S'_{M^{k+1}_{t}}$. The matrix $S_{M^{k}}$ is the two-dimensional concatenation of $\{S'_{M^{k+1}_{t}}\}$ for $t \in \{0,1,2,3\}$. Next, we identify which quadrants $A^{k+1}_{t}$, $B^{k+1}_{u}$, and $C^{k+1}_{v}$ to consider for recursive multiplication. For each quadrant $A^{k+1}_{t} \in \{NW_{A^{k}}, NE_{A^{k}}, SW_{A^{k}}, SE_{A^{k}}\}$ of $A^{k}$, we have identified $S'_{A^{k+1}_{t}}$ (same for for $S'_{B^{k+1}_{u}}$ and $S'_{C^{k+1}_{v}}$, for all $u,v \in \{0,1,2,3\}$). We need to perform the multiplications within $S_{B^{k}}$ and $S_{C^{k}}$ required to calculate $S'_{A^{k+1}_{t}}$. As an example, examine $NW_{A^{k}}$ from Fig. \ref{submatrix_mult_example_mh1}. The sub-matrix $S'_{NW_{A^{k}}}$ is given by the 4-tuple $(A, 5, 3, 3)$. We identify this tuple using Step 2 above. To calculate $S'_{NW_{A^{k}}}=(A,5,3,3)$ of $S_{A^{k}}$, the sub-matrix from $S_{B^{k}}$ given by $(B, 9, 3, 4)$
is to be multiplied by the sub-matrix from $S_{C^{k}}$ given by $(C,
11, 4, 3)$. According to our approach, this will be done in a way so as 
to ensure that the sub-matrices being multiplied at the base case
are \conts, which improves locality and reduces conversion overhead
as described earlier. Because the sub-matrices $(B, 9, 3, 4)$ and
$(C, 11, 4, 3)$ of $S_{B^{k}}$ and $S_{C^{k}}$ touch on all four
quadrants of $B^{k}$ and $C^{k}$, and we want to calculate
$S'_{NW_{A^{k}}}$, all quadrants of $B^{k}$ are to be considered for
multiplication with all quadrants of $C^{k}$. Those are given by the
following sixteen combinations of quadrants from $B^{k}$ and
$C^{k}$:
\begin{equation}
\label{eqn:combs}
\begin{array}{ccc}
NW_{B^{k}} $ and $ NW_{C^{k}} & \quad\quad & SW_{B^{k}} $ and $ NW_{C^{k}}\\
NW_{B^{k}} $ and $ NE_{C^{k}} & \quad\quad & SW_{B^{k}} $ and $ NE_{C^{k}}\\
NW_{B^{k}} $ and $ SW_{C^{k}} & \quad\quad & SW_{B^{k}} $ and $ SW_{C^{k}}\\
NW_{B^{k}} $ and $ SE_{C^{k}} & \quad\quad & SW_{B^{k}} $ and $ SE_{C^{k}}\\
           $     $            & \quad\quad &            $     $           \\
NE_{B^{k}} $ and $ NW_{C^{k}} & \quad\quad & SE_{B^{k}} $ and $ NW_{C^{k}}\\
NE_{B^{k}} $ and $ NE_{C^{k}} & \quad\quad & SE_{B^{k}} $ and $ NE_{C^{k}}\\
NE_{B^{k}} $ and $ SW_{C^{k}} & \quad\quad & SE_{B^{k}} $ and $ SW_{C^{k}}\\
NE_{B^{k}} $ and $ SE_{C^{k}} & \quad\quad & SE_{B^{k}} $ and $ SE_{C^{k}}\\
\end{array}
\end{equation}
All of these are needed to calculate $S'_{NW_{A^{k}}}$. But, to calculate
the sub-matrix $S_{A^{k}}$, we need to find $S'_{NE_{A^{k}}}$,
$S'_{SW_{A^{k}}}$, and $S'_{SE_{A^{k}}}$ in addition to
$S'_{NW_{A^{k}}}$ because $S_{A^{k}}$ is the two-dimensional
concatenation of $\{S_{A^{k+1}_{t}}\}$ for $t \in \{0,1,2,3\}$. Similarly as above, to determine each of
these quadrants of $A^{k}$ requires sixteen combinations of
quadrants from $B^{k}$ and $C^{k}$. In total, to find $S_{A^{k}}$,
we would need up to sixty-four combinations of quadrants from $A^{k}$, $B^{k}$,
and $C^{k}$.

\noindent {\bf Step 4}: For each combination, if $S'_{A^{1}_{t'}}$, $S'_{B^{1}_{u'}}$, and
$S'_{C^{1}_{v'}}$ are not \compat, we extract
\emph{\textbf{compatible}} parts of these sub-matrices and we label them as $S_{A^{1}_{t'}}$,
$S_{B^{1}_{u'}}$, and $S_{C^{1}_{v'}}$ on which the multiplication
proceeds recursively. How to extract compatible parts is beyond the scope of the present manuscript and is left for future work \footnote{We also note that omiting this part of the algorithm does not deflect from its main rationale.}. For now, we concede that omiting it does not divert from the general understanding of the overall algorithm, and that the work requirements for this step can be embedded in that required to perform Steps 1 --> 3 above.

\begin{prop}
If using auxiliary space to peform the matrix additions, and assuming the matrix is of dimensions at most $2^{\alpha} \times 2^{\alpha}$, where $\alpha$ is the machine word-size, the cache oblivious MM using Morton-hybrid order requires asymptotically the same work and critical path lenth as default MM. 
\end{prop}
{\it Proof}:
{\it On work}: The cache-oblivious algorithm is a divide and conquer algorithm. The divide phase introduces two new functions over the default MM algorithm consisting of Steps 1 and 2 above. Each of these steps requires a constant number of arithmetic operations and calls to encoding and extraction procedures. From Sec. 3.5 of \cite{AAArx16a}, we know that each encoding or extraction procedure incurs a constant number of operations assuming the matrix is of dimensions at most $2^{\alpha} \times 2^{\alpha}$, where $\alpha$ is the machine word-size. For the typical value $\alpha = 64$, such matrix sizes are sufficiently large for many applications. It follows that the work of the cache-oblivious algorithm is asymptotically the work of the default algorithm given by $\Theta(n^3)$. The conquer part creates non-overlapping sub-problems in Steps 3 and 4 above whose union yields the original matrix to be multipled.

{\it On parallelism}: All of the extra 64 recursive calls are independent and thus can be cast in parallel. If auxiliary space is available to perform the matrix additions required for each MM, one can also perform addition in parallel using the standard algorithm (Ch. 27 of \cite{CLRS}). Hence, the critical path length of the cache-oblivious algorithm remains that of the default multithreaded algorithm and is known to be $\Theta(\lg^2n)$. 

{\bf Remarks on implications for Parallel Performance:} The sub-matrices at the base case of the recursion are \conts within a row-major sub-matrix, thanks to enforcing \alineds sub-matrices for the recursive division. The Morton-hybrid, cache-oblvious version demonstrates superior performance over the default algorithm, and eliminates the need for Morton-hybrid index conversion when accessing each element in the sub-matrix at the base case, as it can proceed instead with row-major encoding. The implications for parallel performance can be captured using the results from \cite{BlellAl.10}, which reveal that nested parallel algorithms for which the natural sequential execution has low cache complexity will also attain good cache complexity on parallel machines with private or shared caches. In this framework, our adaptation combines improved temporal locality using the Morton-hybrid order for the serial algorithm as well as optimal work and critical path length for the multithreaded version.

\noindent{\bf Performance Analysis} We now verify that the cost of increased recursive MM calls for the cache-oblivious sub-matrix multiplication is significantly compensated for by the improvement in temproal locality thanks to the Morton-hybrid order. We use a Pentium IV of 2.8 GHz processor speed, with an 8 KB L1 cache and a 512 KB L2 cache. It runs linux version 2.6.11 and \texttt{gcc} compiler version 4.0.0. We generate random Morton-hybrid matrices and multiply random sub-matrices of these matrices using both the default and cache-oblivious algorithms. To neutralise the effect of modular aritmetic over finite fields and to be able to exclusively account for the gains induced by the Morton-hybrid order, the random matrices we generate are taken over the binary field. According to \cite{Wise:7At1}, $T' = 32$ is the typical value for the truncation size for block recursive matrix algorithms of floating point entries that shows improvements in cache misses and cycles for Morton-hybrid, default MM. Recall the multiplication of rectangular sub-matrices $S_{A} = S_{B} \cdot S_{C}$, where $A$, $B$ and $C$ are square and in Morton-hybrid order. The dimensions of the square matrices are of no significance, since the multiplication kernel is operating on the rectangular sub-matrices. We thus partition Morton-hybrid matrices of dimensions $N=2048$ and multiply sub-matrices of these Morton-hybrid matrices of varying sizes. Each experiment is distinguished using varying indices $\sigma_{S_{M}}$ of the starting entries of each $S_{M}$ and varying dimensions $r_{S_{M}}$ and $c_{S_{M}}$. 
Because of the variation in sizes across each experiment we do not report on the run-times of each but rather choose to report on the percentage of increase, or decrease, in the number of base case calls made by the cache-oblivious over the default algorithm and the associated percentage of improvement. We record the number of recursive MM calls made to the base case of each of the two algorithms and the total time taken by the overall multiplication to finish. The results are presented in Table \ref{tbl:perc_rec_overhead}. We interpret it using the fifth row, say, as an arbitrary example. Of all 468 experiments run in total, about 9\% of them exhibited about 34\% increase in recursive calls made by the cache-oblivious over the default algorithm. The average, maximum, and minimum percentages of improvement in run-time across this batch of experiments is shown thereafter, and are all staggeringly high. Examining all rows, one can see that no matter what the increase in MM recursive calls has been, this hardly affects the high percentages of improvement. The reductions in cache misses as a result of the cache-oblivious algorithm overwhelm the cost to handle extra recursive calls. 
\begin{table}
  \caption{Percentage Improvement in Runtime of the cache-oblivious algorithm}
  \label{tbl:perc_rec_overhead}
\begin{tabular}{|c|c|c|c|c|}
\hline \textbf{\% Inc. in Calls}  & \textbf{\% of Exp.} &  \textbf{Avg. Imp.} & \textbf{Min. Imp.} & \textbf{Max. Imp.} \\
\hline \textbf{0}           & 6.8  &    96.15  &    95.40   &   97.73   \\
\hline \textbf{3.13}        & 1.7  &    95.94  &    95.40   &   96.59   \\
\hline \textbf{6.25}        & 1.7  &    96.01  &    95.40   &   96.59   \\
\hline \textbf{24.14}       & 1.3  &    95.69  &    95.19   &   96.19   \\
\hline \textbf{34.38}       & 9.4  &    95.93  &    94.83   &   96.61   \\
\hline \textbf{37.5}        & 20.3 &    95.82  &    94.83   &   96.61   \\
\hline \textbf{38.57}       & 0.9  &    95.94  &    95.73   &   96.15   \\
\hline \textbf{41.8}        & 0.9  &    95.75  &    95.74   &   95.76   \\
\hline \textbf{42.77}       & 0.9  &    95.27  &    94.83   &   95.73   \\
\hline \textbf{46.09}       & 0.9  &    95.73  &    94.83   &   96.61   \\
\hline \textbf{80.57}       & 2.6  &    95.88  &    95.48   &   96.15   \\
\hline \textbf{84.77}       & 5.1  &    95.64  &    94.87   &   96.20   \\
\hline \textbf{100}         & 18.6 &    96.03  &    95.35   &   97.73   \\
\hline \textbf{106.25}      & 1.7  &    95.67  &    95.40   &   96.01   \\
\hline \textbf{112.5}       & 1.7  &    95.93  &    95.35   &   96.59   \\
\hline \textbf{168.75}      & 10.3 &    95.60  &    94.83   &   96.61   \\
\hline \textbf{175}         & 5.1  &    95.93  &    94.92   &   96.61   \\
\hline \textbf{300}         & 10.3 &    95.69  &    95.35   &   97.73   \\
\hline
\end{tabular}\end{table}

\end{document}